\documentclass[12pt, a4wide]{amsart}
\usepackage{mathrsfs}
\usepackage{amsmath,amssymb}
\newcommand{\be}{\begin{equation}}
\newcommand{\ee}{\end{equation}}
\newcommand{\bal}{\begin{aligned}}
\newcommand{\eal}{\end{aligned}}
\newcommand{\bee}{\begin{equation*}}
\newcommand{\eee}{\end{equation*}}
\usepackage{color}

\def\cwedge{\bigcirc\kern-1.07em\wedge\ }

\newtheorem{theorem}{Theorem}[section]

\newtheorem{proposition}[theorem]{Proposition}
\newtheorem{corollary}[theorem]{Corollary}

\numberwithin{equation}{section}

\begin{document}

\title[{\it curvature identity}]{A curvature identity on a 6-dimensional Riemannian Manifold and its applications}

\author[Y. Euh]{Yunhee Euh}

\address{National Institute of Mathematical Sciences,
Deajeon 305-811, Korea}
\email{yheuh@nims.re.kr}

\author[J. H. Park]{JeongHyeong Park}
\address{Department of Mathematics,
    Sungkyunkwan University,
    Suwon 440-746, Korea}
\email{parkj@skku.edu}
\author[K. Sekigawa]{Kouei Sekigawa}
\address{ Department of Mathematics,
    Niigata University,
    Niigata 950-2181, Japan}
\email{sekigawa@math.sc.niigata-u.ac.jp}


\begin{abstract}

We derive a curvature identity that holds on any 6-dimensional
Riemannian manifold, from the Chern-Gauss-Bonnet theorem for
a 6-dimensional closed Riemannian manifold. We also introduce some
applications of this curvature identity.
\end{abstract}

\maketitle
\section{Introduction}

M. Berger \cite{Ber} derived a curvature identity on 4-dimensional
compact Riemannian manifolds from the Chern-Gauss-Bonnet theorem
based on the well-known fact that the Euler number is a topological
invariant. We demonstrated that the obtained
curvature identity holds on any 4-dimensional Riemannian manifold
which is not necessarily compact \cite{EPS2}. Further, Gilkey, Park
and Sekigawa extended the result to the higher dimensional setting,
the pseudo-Riemannian setting, manifolds with boundary setting and
the K\"ahler setting \cite{GPS1, GPS2, GPS3, GPS4}. In this paper,
we shall give a curvature identity explicitly which holds on any
6-dimensional Riemannian manifold using methods similar to those
used in the 4-dimensional Chern-Gauss-Bonnet theorem and also
provide some applications of the obtained curvature identity. More
precisely, we derive a symmetric 2-tensor valued curvature identity of
degree 6 which holds on any 5-dimensional Riemannian manifold, from
which a scalar-valued curvature identity can be derived
(\cite{GPS1}, Lemma 1.2 (3)). Furthermore, we derive a symmetric
2-tensor valued curvature identity of degree 6 on 4-dimensional
Riemannian manifolds from the curvature identity on 5-dimensional
Riemannian manifolds. Based on these obtained identities, we shall
also discuss a question that arose in \cite{CGW} related to the
Lichnerowicz conjecture for a harmonic manifold ``a harmonic
manifold is locally symmetric". The original Lichnerowicz conjecture
is the one for 4-dimensional harmonic manifolds which was proved by
Walker (\cite{W} and Corollary 1.2 in \cite{CGW}). The
Lichnerowicz conjecture was refined by Ledger since he showed that a
locally symmetric manifold is harmonic if and only if it is locally
isometric to a Euclidean space or a rank one symmetric space
\cite{Le2}. Concerning the Lichnerowicz conjecture, Szab\'o
\cite{Sz} proved that the conjecture is true on the compact harmonic
manifolds. For the non-compact case, Damek and Ricci \cite{BTV, DR} provided
the counterexample demonstrating that the Lichnerowicz conjecture
is not true for case dimensions $\geqq 7$. As mentioned above,
the Lichnerowicz conjecture is true for the 4-dimensional case.
Further, Nikolayevsky \cite{Ni} showed that the Lichnerowicz
conjecture refined by Ledger is also true for the
5-dimensional case. Presently, to the best of our knowledge, the
Lichnerowicz conjecture is still open for the 6-dimensional case. In
the present paper, we provide an another proof of the Lichnerowicz
conjecture the refined version by Ledger for the 4-dimensional case and
a brief review of the proof of the Lichnerowicz conjecture for
5-dimensional case by Nikolayevsky
under slightly general settings. For more detailed information
concerning the Lichnerowicz conjecture, refer to \cite{Li, Ni, Pa,
Pe}.

\section{Preliminaries}

In this section, we shall prepare several fundamental concepts,
terminologies and notational conventions. In the present paper, we
shall  adopt similar notational conventions as those used in
\cite{GPS1}. We denote by $\mathcal{I}_{m,n}$ the space of scalar
invariant local formulas and by $\mathcal{I}^2_{m,n}$  the space of
symmetric 2-tensor valued invariant local formulas, respectively,
defined in the category of all Riemannian manifolds of dimension $m$
and of degree $n$. We note that $\mathcal{I}_{m,n}=\{0\}$ and
$\mathcal{I}^2_{m,n}=\{0\}$ if $n$ is odd. We denote by $r$ the
restriction map $r:\mathcal{I}_{m,n}\rightarrow \mathcal{I}_{m-1,n}$
(resp. $r:\mathcal{I}^2_{m,n}\rightarrow \mathcal{I}^2_{m-1,n}$)
given by restricting the summation to range from 1 to $m-1$.

Now, let $M=(M,g)$ be an $m$-dimensional Riemannian manifold and $\nabla$ be
the Levi-Civita connection of $g$. We assume that the curvature tensor $R$ is defined by
\begin{equation}\label{Def:curv}
R(X,Y)Z=[\nabla_X,\nabla_Y]Z-\nabla_{[X,Y]}Z
\end{equation}
for $X$, $Y\in \mathfrak{X}(M)$, where $\mathfrak{X}(M)$ denotes the
Lie algebra of all smooth vector fields on $M$. We also denote the
Ricci tensor and the scalar curvature of $M$ by $\rho$ and $\tau$,
respectively.  Let $\{e_i\}=\{e_1, e_2, \cdots e_m\}$ be a local
orthonormal frame and $\{e^i\}$ be a dual frame field. Throughout
the present paper, we assume that the components of the tensor
fields are relative to a local orthonormal frame $\{e_i\}$
and also adopt the Einstein convention on sum over repeated indices
unless otherwise specified. Further, we denote by $R_{abcd;i}$,
$R_{abcd;ij},\cdots$ the components of the covariant derivatives
of the curvature tensor $R=(R_{abcd})$ with respect to the
Levi-Civita connection $\nabla$. The following theorems play
fundamental roles in our forthcoming discussion.

\begin{theorem}\label{thm-2.1}{\rm\cite{GPS1}}

\noindent (1) $r:\mathcal{I}_{m,n}\rightarrow\mathcal{I}_{m-1,n}$ is surjective.

\smallbreak\noindent (2) If $n$ is even and if $m>n$, then $r:\mathcal{I}_{m,n}\rightarrow\mathcal{I}_{m-1,n}$  is bijective.

\smallbreak\noindent (3) Let $m$ be even. Then $ker\{r:\mathcal{I}_{m,m}\rightarrow\mathcal{I}_{m-1,m}\}=E_{m,m}\cdot\mathbb{R}$, where $E_{m,n}\in \mathcal{I}_{m,n}$ is the Pfaffian form defined by
\begin{equation}
E_{m,n}:=\sum_{i_1,...,i_n,j_1,...,j_n=1}^m{{R_{i_1i_2j_2j_1}...R_{i_{n-1}i_{n}j_nj_{n-1}}}}g(e^{i_1}\wedge...\wedge
e^{i_n},e^{j_1}\wedge...\wedge e^{j_n}),
\end{equation}
\end{theorem}

\begin{theorem}\label{thm-2.2}{\rm\cite{GPS1}}

\noindent(1) $r:\mathcal{I}_{m,n}^2\rightarrow\mathcal{I}_{m-1,n}^2$ is surjective.

\smallbreak\noindent(2) If $n$ is even and if $m>n+1$, then $r:\mathcal{I}_{m,n}^2\rightarrow\mathcal{I}_{m-1,n}^2$ is bijective.

\smallbreak\noindent(3)  If $m$ is even, then
$ker\{r:\mathcal{I}_{m+1,m}^2\rightarrow\mathcal{I}_{m,m}^2\}=
T_{m,m}^2\cdot\mathbb{R}$, where $T_{m,n}^2\in \mathcal{I}^2_{m,n}$
is the Pfaffian defined by
\begin{equation}
\begin{aligned}
T_{m,n}^2 :=\sum_{i_1,...,i_{n+1},j_1,...,j_{n+1}=1}^mR_{i_1i_2j_2j_1}...R_{i_{n-1}i_{n}j_{n}j_{n-1}}e^{i_{n+1}}\circ
e^{j_{n+1}}\\
\qquad\qquad\qquad\qquad\times
g(e^{i_1}\wedge...\wedge e^{i_{n+1}},e^{j_1}\wedge...\wedge e^{j_{n+1}})\,.
\end{aligned}
\end{equation}

\end{theorem}

\section{The universal curvature identity}

Let $M=(M,g)$ be a 6-dimensional compact oriented Riemannian
manifold. Then, it is well-known that the Euler number of $M$ is
given by the following integral formula, namely, the
Chern-Gauss-Bonnet theorem.
\begin{theorem}\label{thm:3.1}{\rm\cite{Sa}}
    \be\label{eq:G_B_6}
    \bal
    \chi(M)&=\frac{1}{2^9\pi^3 3!}\int_ME_{6,6}dv_g\\
    &=\frac{1}{384\pi^3}\int_M \{\tau^3-12\tau|\rho|^2+3\tau|R|^2+16\rho_{ab}\rho_{ac}\rho_{bc}-24\rho_{ab}\rho_{cd}R_{acbd}\\
    &-24\rho_{uv}R_{uabc}R_{vabc}+8R_{abcd}R_{aucv}R_{bvdu}-2R_{abcd}R_{abuv}R_{cduv}\}dv_g.
    \eal
    \ee
\end{theorem}
We here set
$$\hat{R}\equiv R_{abcd}R_{abuv}R_{cduv} $$ and $$\mathring{R}\equiv R_{abcd}R_{aucv}R_{budv}.$$
We note that identity \eqref{eq:G_B_6} is rearranged by our
setting: the curvature of \cite{Sa} has a negative sign
difference to ours and the term $8R_{abcd}R_{aucv}R_{bvdu}$
has been changed by using the first Bianchi identity to
$8\mathring{R}-2\hat{R}$.

Now, we regard the right hand side of \eqref{eq:G_B_6} as a
functional $\mathcal{F}$ on the space $\mathfrak{M}(M)$ of all
Riemannian metrics on $M$. Let $h$ be any symmetric
$(0,2)$-tensor field in $M$ and consider a one-parameter deformation
of $ g$ by $g(t)=g+th$ for any $g\in\mathfrak{M}(M)$.  Since the
Euler number $\chi(M)$ is a topological invariant of $M$,
$\mathcal{F}$ does not depend on the choice of Riemannian metrics on
$M$, so we have
\begin{equation}\label{eq:3.2}
    0=\frac{d}{dt}\Big|_{t=0}\mathcal{F}{(g(t))}=0.
\end{equation}
This holds for any symmetric (0,2)-tensor field $h$ on $M$.
Applying the similar arguments as in \cite{EPS2}, taking account
of \eqref{eq:G_B_6} and \eqref{eq:3.2}, we have the
following equality as the corresponding Euler-Lagrange equation for
the functional $\mathcal{F}$:

 \begin{equation}\label{eq:main eq}
    \begin{aligned}
    &\frac{1}{2}({\tau^3}+{3}\tau|R|^2-12\tau|\rho|^2+16\rho_{ab}\rho_{bc}\rho_{ca}-24\rho_{ab}\rho_{cd}R_{acbd}-24\rho_{uv}R_{abcu}R_{abcv}+8\mathring{R}- 4\hat{R})g_{ij}\\
    &-3\tau^2\rho_{ij}-3|R|^2\rho_{ij}+12|\rho|^2\rho_{ij}+12\tau\rho_{ia}\rho_{ja}+12\tau\rho_{ab}R_{iabj}-6\tau R_{iabc}R_{jabc}\\
   &-24\rho_{ia}\rho_{jb}\rho_{ab}-24\rho_{ac}\rho_{bc}R_{iabj}+24\rho_{aj}\rho_{cd}R_{acid}+24\rho_{ai}\rho_{cd}R_{acjd}+24\rho_{ab}R_{icdj}R_{acbd}  \\
    &+48\rho_{cd}R_{iabc}R_{jabd}+6\rho_{jd}R_{abci}R_{abcd} +6\rho_{id}R_{abcj}R_{abcd}+12 \check{R}_{ij} +12 \hat{R}_{ij} -24
    \mathring{R}_{ij}=0,
    \end{aligned}
    \end{equation}
where $\check{R}_{ij}=R_{iuvj}R_{abcu}R_{abcv}$, $\hat{R}_{ij}=
R_{ibac}R_{jbuv}R_{acuv}$ and
$\mathring{R}_{ij}=R_{iabc}R_{jubv}R_{aucv}$. We here omit the
detailed calculation. From \eqref{eq:main eq} and Theorem
\ref{thm:3.1}, taking account of the results of (\cite{EGPS},
Theorem 1.2) and (\cite{EPS2}, Main theorem), we have the
following.

\begin{theorem}
The curvature identity \eqref{eq:main eq} holds on any 6-dimensional Riemannian manifold $M=(M,g)$ which is not necessarily compact and further it is universal in $\mathcal{I}_{6,6}^{2}$.
\end{theorem}
\noindent Especially, we have the following.
\begin{corollary}\label{cor:3.3}
Let $M=(M,g)$ be a 6-dimensional Einstein manifold. Then the following identity holds on $M$:
\begin{equation}\label{eq:Cor}
(-\tau|R|^2+4\mathring{R}{-2\hat{R}})g_{ij}+12\check{R}_{ij}+12\hat{R}_{ij}-24\mathring{R}_{ij}+4\tau
R_{iabc}R_{jabc}=0.
\end{equation}
\end{corollary}
\vskip0.15cm We note that the curvature identity \eqref{eq:main eq}
can also be obtained by making use of the equality $T_{6,6}^{2}= 0$
from Theorem 2.2 (3). However, we derived the same identity
\eqref{eq:main eq} without adopting this method in this
paper. Further, we note that the curvature identity is universal in
the same form for any 6-dimensional pseudo-Riemannian manifold
\cite{GPS2}.

\section{Derived curvature identities on 4- and 5-dimensional \\Riemannian manifolds}
In this section, we shall provide further curvature identities on
4- and 5-dimensional Riemannian manifolds derived from the
curvature identity \eqref{eq:main eq} on 6-dimensional Riemannian
manifolds.

Now, let $M=(M,g)$ be a 5-dimensional Riemannian manifold and
$\bar{M}=(M\times \mathbb{R}, g \oplus 1)$ be the Riemannian product
of $M=(M,g)$ and a real line $\mathbb{R}$. Then, applying Theorem
\ref{thm:3.1} to the Riemannian manifold $\bar{M}=(M\times \mathbb{R}, g \oplus 1)$, we see that the following curvature identity
\begin{equation}\label{eq:4.1}
\begin{aligned}
    &\tau^3-12\tau|\rho|^2+3\tau|R|^2+16\rho_{ab}\rho_{bc}\rho_{ca}\\
    &-24\rho_{ab}\rho_{cd}R_{acbd}-24\rho_{uv}R_{abcu}R_{abcv}+8\mathring{R}-4\hat{R}=0
\end{aligned}
\end{equation}
holds on $M$ and further, it is universal in $\mathcal{I}_{5,6}$ (\cite{GPS1},  Lemma 1.2 (3)).

On one hand, taking account of Theorem \ref{thm-2.2} (1), we see
that \eqref{eq:main eq} holds on $M$ in the same form by restricting
the range of the indices from 1 to 5. Therefore, from \eqref{eq:main
eq} and \eqref{eq:4.1}, we have the following.
\begin{theorem}\label{thm:4.2}
Let $M=(M,g)$ be a 5-dimensional Riemannian manifold. Then, in addition to \eqref{eq:4.1}, the following identity
  \begin{equation}\label{eq:4.2}
\begin{aligned}
  &\tau^2\rho_{ij}+|R|^2\rho_{ij}-4|\rho|^2\rho_{ij}-4\tau\rho_{ia}\rho_{ja}-4\tau\rho_{ab}R_{iabj}+2\tau R_{iabc}R_{jabc}\\
   &+8\rho_{ia}\rho_{jb}\rho_{ab}+8\rho_{ac}\rho_{bc}R_{iabj}-8\rho_{aj}\rho_{cd}R_{acid}-8\rho_{ai}\rho_{cd}R_{acjd}-8\rho_{ab}R_{icdj}R_{acbd} \\
    &-16\rho_{cd}R_{iabc}R_{jabd} -2\rho_{jd}R_{abci}R_{abcd}-2\rho_{id}R_{abcj}R_{abcd}-4 \check{R}_{ij} -4\hat{R}_{ij}+8\mathring{R}_{ij}=0
\end{aligned}
\end{equation}
holds on M.
\end{theorem}
\noindent{\bf Remark 1} Tranvecting \eqref{eq:4.2} with $g_{ij}$, we
may also obtain \eqref{eq:4.1}.\\
From Theorem \ref{thm:4.2}, we have the following.
\begin{corollary}\label{cor:4.2}
Let $M=(M,g)$ be a 5-dimensional Einstein manifold. Then, we have
\begin{equation}\label{eq:5-dim Eins}
\Big(\frac{\tau^3}{25} +\frac{\tau}{5}|R|^2\Big)g_{ij}-2\tau
R_{iabc}R_{jabc}-4\check{R}_{ij}-4\hat{R}_{ij}+8\mathring{R}_{ij}=0.
\end{equation}
\end{corollary}
\noindent From \eqref{eq:4.2}, taking account of Theorem
\ref{thm-2.2} (1) and Equation (1.2) in \cite{EPS3}, we have
the following.
\begin{corollary}\label{cor:4.3}
Let $M=(M,g)$ be a  4-dimensional Riemannian manifold. Then, the identity \eqref{eq:4.2} holds in the same form by restricting
the range of the indices from 1 to 4 and further, it is universal in $\mathcal{I}^2_{4,6}$. Especially, if $M$ is Einstein, the identity reduces to the following identity:
\begin{equation*}
\Big(\frac{\tau^3}{8}-\frac{3}{4}\tau|R|^2\Big)g_{ij}-4\hat{R}_{ij}+8\mathring{R}_{ij}=0.
\end{equation*}
\end{corollary}
\noindent {{Here, we shall call a 6-dimensional, 5-dimensional
    and 4-dimensional Riemannian manifold satisfying
    the curvature identities in the Corollaries 3.3, 4.2
    and 4.3, a 6-dimensional, 5-dimensional and
    4-dimensional weakly Einstein manifold of degree 6,
    respectively.
  Based on our current work, the definition 4-dimensional weakly
    Einstein manifold introduced in our paper \cite{EPS2, EPS3}
    may be made more precise and the definition becomes
    4-dimensional weakly Einstein manifold
    of degree 4. We note that Arias-Marco and Kowalski recently obtained a
          classification theorem for 4-dimensional homogeneous
          weakly Einstein manifolds \cite{AK}}}.

\vspace{3mm}

\section{A generalization of harmonic manifolds}
An $m$-dimensional Riemannian manifold $M=(M,g)$ is called a {\it locally
harmonic} manifold (briefly, harmonic manifold) if, for every point
$p\in M$, the volume density function $\theta_p(q)=\sqrt{\det
(g_{ij})(q)}$ is a radial function in a normal neighborhood
$U_p=U_p(x^1,\cdots,x^m)$ centered at $p$, where
$g_{ij}=g(\partial/\partial x_i,\partial/\partial x_j)$, namely,
there exists a positive real number $\varepsilon(p)$ and a smooth
function $\Theta_p:[0,\varepsilon(p))\rightarrow M$ such that
$\theta_p(q)=\Theta_p(d(p,q))$ for $q\in U_p$ where $d(p,q)$ is a
distance from $p$ to $q$. We here note that there are several
equivalent definitions for harmonic manifolds \cite{Be}. A locally
Euclidean space and a locally rank one symmetric space are harmonic
manifolds. Concerning the converse, there is a well-known conjecture
known as the Lichnerowicz conjecture that every harmonic manifold is
locally isometric to a Euclidean space or a rank one symmetric
space. Copson and Ruse \cite{CR}, Lichnerowicz \cite{Li44} and
Ledger \cite{Le} have shown that each harmonic manifold must satisfy
an infinite sequence $\{H_n\}_{n=1,2,\cdots}$ of conditions on the
curvature tensor and its covariant derivatives. The first three of these conditions are
given as follows \cite{Be,Wa}:

\begin{equation*}
\begin{aligned}
H_1:&\quad R_{aija}=\Lambda_1 g_{ij},\\
H_2:&\quad\mathfrak{S}(R_{aijb}R_{bkla})=\Lambda_2\mathfrak{S}(g_{ij}g_{kl}),\\
H_3:&\quad\mathfrak{S}(32R_{aijb}R_{bklc}R_{cuva}-9R_{aijb;k}R_{buva;l})=\Lambda_3\mathfrak{S}(g_{ij}g_{kl}g_{uv}),
\end{aligned}
\end{equation*}
where each $\Lambda_n$ $(n=1,2,3)$ is a constant and $\mathfrak{S}$
denotes the summation taken over all permutations of the free
indices appearing inside the parenthesis. From the condition $H_1$,
it follows immediately that a harmonic manifold is Einstein and
hence real analytic as a Riemannian manifold. We may note that the
conditions $H_{1},H_{2},H_{3}$ are equivalent to the following
conditions ${H'_{1}},{H'_{2}},{H'_3}$ respectively \cite{CGW}:

\begin{equation*}
\begin{aligned}
{H'_1}: &\quad R_{axxa} = \Lambda_1 |x|^2, \\
{H'_2}: &\quad R_{axxb}R_{bxxa} = \Lambda_2 |x|^4, \\
{H'_3}: &\quad 32R_{axxb}R_{bxxc}R_{cxxa} - 9R_{axxb;x}R_{bxxa;x} = \Lambda_3 |x|^6,
\end{aligned}
\end{equation*}
for any $x = {\xi}^{i}e_{i} \in T_{p}{M}$ at $p \in M$, where
$R_{axxb} = R_{aijb}{\xi}^{i}{\xi}^{j}$ and $R_{axxb;c} =
R_{aijb;k}{\xi}^{i}{\xi}^{j}{\xi}^{k}$.
 \vskip0.15cm
\noindent{\bf Remark 2} The condition $H_3$ in \cite{CGW} is
incorrect (\cite{Be}, pp.162) and should be changed for the above $H'_{3}$.
 \vskip0.15cm

 In \cite{CGW},
Carpenter, Gray and Willmore raised the following question:
\vskip0.2cm
\noindent{\bf Question A.} Does a Riemannian manifold $M=(M,g)$ exist which satisfies some but
not all of the conditions $\{H_n\}_{n=1,2,\cdots}$?

\vskip0.2cm
\noindent Concerning Question A, they discussed the case where
$M=(M,g)$ is a non-flat locally symmetric space satisfying the
condition $H_1$ and some other condition $H_k$ and obtained some
partial answers to the question (\cite{CGW}, Theorem 1.1). Taking
account of these observations, it seems worthwhile to consider the
Question A under a more general setting.

Now, we shall define a generalization of harmonic manifolds.
 \vskip0.2cm
\noindent {\bf Definition 5.1} A Riemannian manifold $M = (M,g)$
satisfying the conditions $\{H_n\}_{n=1,\cdots,k}$ is called an
{\it{asymptotic harmonic manifold up to order $k$.}}
 \vskip0.2cm
By the above definition, it follows immediately that an asymptotic
harmonic manifold up to order $k$ is an asymptotic harmonic manifold
up to order for any $\ell(1 \leq\ell < k)$. Further, we may
check that a locally symmetric asymptotic harmonic manifold up to
order $k$ is $k$-stein \cite{CGW, Le, Le2}.

Let $M=(M,g)$ be an $m$-dimensional asymptotic harmonic manifold
up to order 2. Then, we have
\begin{equation}\label{eq:new4.1}
\rho_{ij}=\Lambda_1g_{ij}\quad\text{ and hence,
}\Lambda_1=\frac{\tau}{m},
\end{equation}
and
\begin{equation}\label{eq:new4.2}
\begin{aligned}
&R_{aijb}R_{aklb}+R_{aijb}R_{alkb}+R_{aikb}R_{ajlb}+R_{aikb}R_{aljb}+R_{ailb}R_{akjb}+R_{ailb}R_{ajkb}\\
&=2\Lambda_2(g_{ij}g_{kl}+g_{ik}g_{jl}+g_{jk}g_{il}).
\end{aligned}
\end{equation}
Transvecting \eqref{eq:new4.2} with $g_{kl}$ and taking account of \eqref{eq:new4.1}, we have
\begin{equation}\label{eq:new4.3}
2\Big(\frac{\tau}{m}\Big)^2g_{ij}+3R_{iabc}R_{jabc}=2(m+2)\Lambda_2g_{ij}.
\end{equation}
From \eqref{eq:new4.3}, we have
\begin{equation}\label{eq:new4.4}
\Lambda_2=\frac{1}{2m(m+2)}\Big(\frac{2\tau^2}{m}+3|R|^2\Big).
\end{equation}
Thus, from \eqref{eq:new4.4}, it follows immediately that $|R|^2$ is constant on $M$. Further, from \eqref{eq:new4.3} and \eqref{eq:new4.4}, we have
\begin{equation}\label{eq:new4.5}
\begin{aligned}
R_{iabc}R_{jabc}&=\frac{1}{3}\Big\{\frac{1}{m}\Big(\frac{2\tau^2}{m}+3|R|^2\Big)-\frac{2\tau^2}{m^2}\Big\}g_{ij}\\
&=\frac{1}{m}|R|^2g_{ij},
\end{aligned}
\end{equation}
and hence, $M$ is a super-Einstein manifold with constant $|R|^2$
(\cite{BV, GW}).\vskip0.1cm
 \vskip0.2cm

\noindent{\bf Remark 3} By definition, an
$m(\geqq3)$-dimensional asymptotic harmonic manifold up to order
2 is a 2-stein manifold with constant $|R|^2$
\cite{CGW}. It is known that for each 2-stein manifold of
dimension $m(\ne4)$, $|R|^2$ is constant. An explicit example of
4-dimensional 2-stein manifold with non-constant $|R|^2$ has been
provided in \cite{CPS}. It is also known that every 2-stein manifold
is super-Einstein. We may reconfirm this fact by the above equality
\eqref{eq:new4.5}.
 \vskip0.2cm

From \eqref{eq:new4.2}, taking account of \eqref{eq:new4.4}, we
may show the following.
\begin{proposition}\label{prop4.1}
Let $M = (M,g)$ be an $m(\geqq 3)$-dimensional non-flat
asymptotic harmonic manifold up to order 2. Then, $M$ is
irreducible.
\end{proposition}

The following identity holds in arbitrary. We shall use it to
derive the Lichnerowicz formula
\begin{equation}\label{eq:new4.12}
\begin{aligned}
&(R_{iabc}R_{jabc})_{;kk}\\
&=2B_{ij}+8\mathring{R}_{ij}+2\hat{R}_{ij} + 4\rho_{cd}R_{iabc}R_{jabd}\\
&+2\rho_{ic;{ab}}R_{jabc}+2\rho_{jc;ab}R_{iabc}+2\rho_{ab;ic}R_{jabc}+2\rho_{ab;jc}R_{iabc}.
\end{aligned}
\end{equation}
Especially, if the Riemannian manifold $M = (M,g)$ is Einstein with constant $|R|^2$, from \eqref{eq:new4.12},
we have easily

\begin{equation}\label{eq:new4.14}
|\nabla R|^2=-4\mathring{R}-\hat{R}-\frac{2\tau}{m}|R|^2.
\end{equation}

In the sequel, we assume that every Riemannian manifold $M=(M,g)$ is
an $m(\geqq 4)$-dimensional asymptotic harmonic manifold up to
order 3 unless otherwise specified. Then, from the condition
$H_3$, taking account of \eqref{eq:new4.5}, we have
\begin{equation}\label{eq:new4.6}
g_{kl}g_{uv}\mathfrak{S}(R_{aijb;k}R_{buva;l})=48(A_{ij}+2B_{ij}),
\end{equation}
\begin{equation}\label{eq:new4.7}
\begin{aligned}
g_{kl}g_{uv}\mathfrak{S}(R_{aijb}R_{bklc}R_{cuva})
&=48\Big(\frac{\tau^3}{m^3}g_{ij}+\frac{3}{2}\check{R}_{ij}-\frac{7}{2}\hat{R}+\mathring{R}_{ij}+\frac{3}{m}\tau R_{iabc}R_{jabc}\Big)\\
&=48\Big\{-\frac{7}{2}\hat{R}_{ij}+\mathring{R}_{ij}+\Big(\frac{\tau^3}{m^3}+\frac{9\tau}{2m^2}|R|^2\Big)g_{ij}\Big\},
\end{aligned}
\end{equation}
\begin{equation}\label{eq:new4.8}
g_{kl}g_{uv}\mathfrak{S}(g_{ij}g_{kl}g_{uv})=48(m+2)(m+4)g_{ij},
\end{equation}

\noindent
where $A_{ij}=R_{abcd;i}R_{abcd;j}$ and $B_{ij}=R_{ibcd;a}R_{jbcd;a}$.
Thus, from $H_3$ and \eqref{eq:new4.6}$\sim$\eqref{eq:new4.8}, we have
\begin{equation}\label{eq:5.9}
\begin{aligned}
&32\Big\{-\frac{7}{2}\hat{R}_{ij}+\mathring{R}_{ij}+\Big(\frac{\tau^3}{m^3}+\frac{9\tau}{2m^2}|R|^2\Big)g_{ij}\Big\}-9(A_{ij}+2B_{ij})\\
&=(m+2)(m+4)\Lambda_3g_{ij}
\end{aligned}
\end{equation}
Multiplying \eqref{eq:5.9} by $m$, we have the following
equation:
\begin{equation}\label{eq:new4.9}
\begin{aligned}
&32\Big\{-\frac{7m}{2}\hat{R}_{ij}+m\mathring{R}_{ij}+\Big(\frac{\tau^3}{m^2}+\frac{9\tau}{2m}|R|^2\Big)g_{ij}\Big\}-9m(A_{ij}+2B_{ij})\\
&=m(m+2)(m+4)\Lambda_3g_{ij}.
\end{aligned}
\end{equation}
Transvecting \eqref{eq:5.9} with $g_{ij}$, we further have
\begin{equation}\label{eq:new4.10}
32\Big(-\frac{7}{2}\hat{R}+\mathring{R}+\frac{\tau^3}{m^2}+\frac{9\tau}{2m}|R|^2\Big)-27|\nabla
R|^2=
 m(m+2)(m+4)\Lambda_3.
\end{equation}
Thus, from \eqref{eq:new4.9} and \eqref{eq:new4.10}, we have
\begin{equation*}
\begin{aligned}
&9m(A_{ij}+2B_{ij})-27|\nabla R|^2g_{ij}\\
&=32\Big\{-\frac{7m}{2}\hat{R}_{ij}+m\mathring{R}_{ij}+\Big(\frac{\tau^3}{m^2}+\frac{9\tau}{2m}|R|^2\Big)g_{ij}\Big\}+27|\nabla R|^2g_{ij}\\
&-32\Big(-\frac{7}{2}\hat{R}+\mathring{R}+\frac{\tau^3}{m^2}+\frac{9\tau}{2m}|R|^2\Big)g_{ij},
\end{aligned}
\end{equation*}
and hence,
\begin{equation}\label{eq:new4.11}
9m(A_{ij}+2B_{ij})-27|\nabla R|^2g_{ij}=32(m\mathring{R}_{ij}-\mathring{R}g_{ij})-112(m\hat{R}_{ij}-\hat{R}g_{ij}).
\end{equation}
\noindent From \eqref{eq:new4.12}, taking account of
\eqref{eq:new4.5}, we have
\begin{equation}\label{eq:new4.13}
B_{ij}=-4\mathring{R}_{ij}-\hat{R}_{ij}-\frac{2\tau}{m^2}|R|^2g_{ij}.
\end{equation}
Thus, from \eqref{eq:new4.1}, \eqref{eq:new4.4}, \eqref{eq:new4.14}
and \eqref{eq:new4.10}, we have the following.
\begin{proposition}\label{prop:4.4}
Let $M=(M,g)$ be an $m$-dimensional asymptotic harmonic manifold
up to order 3. Then, $M$ is a 2-stein manifold with constant
$|R|^2$, and further, $|\nabla R|^2+\hat{R}+4\mathring{R}$,
$27|\nabla R|^2+112\hat{R}-32\mathring{R}$ are constant and
hence, $17\hat{R}-28\mathring{R}$ is constant on $M$.
\end{proposition}

\noindent{\bf Remark 4} Proposition 6.68 in \cite{Be} should be corrected as above. \\
Here, we set
\begin{equation}\label{eq:new4.15}
\begin{aligned}
\alpha_{ij}=&A_{ij}-\frac{1}{m}|\nabla R|^2g_{ij},\\
\beta_{ij}=&B_{ij}-\frac{1}{m}|\nabla R|^2g_{ij},\\
\hat{\gamma}_{ij}=&\hat{R}_{ij}-\frac{1}{m}\hat{R}g_{ij},\\
\mathring{\gamma}_{ij}=&\mathring{R}_{ij}-\frac{1}{m}\mathring{R}g_{ij}.
\end{aligned}
\end{equation}
Then, from  \eqref{eq:new4.11} and \eqref{eq:new4.15}, we have
\begin{equation}\label{eq:new4.16}
9(\alpha_{ij}+2\beta_{ij})=32\mathring{\gamma}_{ij}-112\hat{\gamma}_{ij},
\end{equation}
and hence, from \eqref{eq:new4.13}$\sim$\eqref{eq:new4.15},
we have
\begin{equation*}
\begin{aligned}
B_{ij}&=-4\mathring{R}_{ij}-\hat{R}_{ij}-\frac{2\tau}{m^2}|R|^2g_{ij}\\
&=-4\mathring{\gamma}_{ij}-\hat{\gamma}_{ij}-\frac{1}{m}(4\mathring{R}+\hat{R}+\frac{2\tau}{m}|R|^2)g_{ij}\\
&=-4\mathring{\gamma}_{ij}-\hat{\gamma}_{ij}+\frac{1}{m}|\nabla R|^2
g_{ij}.
\end{aligned}
\end{equation*}
Thus we obtain
\begin{equation}\label{eq:new4.17}
\beta_{ij}=-4\mathring{\gamma}_{ij}-\hat{\gamma}_{ij}.
\end{equation}
Hence, we have the following.
\begin{proposition}\label{prop:4.6}
Let $M=(M,g)$ be an $m$-dimensional asymptotic harmonic manifold
up to order 3. Then the following equalities hold:
\begin{equation*}
9\alpha_{ij}=104\mathring{\gamma}_{ij}-94\hat{\gamma}_{ij},\qquad
\beta_{ij}=-4\mathring{\gamma}_{ij}-\hat{\gamma}_{ij}.
\end{equation*}
\end{proposition}

\subsection {4-dimensional asymptotic harmonic manifolds }

Let $M = (M,g)$ be a 4-dimensional asymptotic harmonic manifold
up to order 3.  Then, since $M$ is a 2-stein manifold (with
constant $|R|^2$) for each point $p\in M$, we may choose a
Singer-Thorpe basis $\{e_i\}=\{e_1,e_2,e_3,e_4\}$ such that
\begin{equation}\label{eq:new4.21}
\begin{gathered}
R_{1212}=R_{3434}= a,\quad R_{1313}=R_{2424} = b,\quad R_{1414}=R_{2323} = c,\\
R_{1234}=\alpha, \quad R_{1342}= \beta, \quad R_{1423} =
\gamma.
\end{gathered}
\end{equation}
satisfying $\alpha+\beta+\gamma=0$ and $\alpha=a+\frac{\tau}{12}$,
$\beta=b+\frac{\tau}{12}$, $\gamma=c+\frac{\tau}{12}$ (or
$-\alpha=a+\frac{\tau}{12}$, $-\beta=b+\frac{\tau}{12}$,
$-\gamma=c+\frac{\tau}{12}$) \cite{SV}.
\newline Without
loss of essentiality, it suffices to consider the case,
\begin{equation}\label{eq:new4.22}
\alpha=a+\frac{\tau}{12},\quad \beta=b+\frac{\tau}{12},\quad \gamma=c+\frac{\tau}{12}.
\end{equation}
Then, by straightforward calculation, we obtain
\begin{equation}\label{eq:new4.23}
\begin{gathered}
\tau=-4(a+b+c),\qquad |R|^2=\frac{5}{6}\tau^2-32(ab+bc+ca),\\
\hat{R}_{ij}=\frac{1}{4}\hat{R}g_{ij}\quad(\hat{\gamma}_{ij}=0),
\qquad \hat{R}=192abc+32\tau(ab+bc+ca)-\frac{7}{12}\tau^3,\\
\mathring{R}_{ij}=\frac{1}{4}\mathring{R}g_{ij}\quad(\mathring{\gamma}_{ij}=0),\qquad\mathring{R}=96abc+4\tau(ab+bc+ca)-\frac{\tau^3}{24}.
\end{gathered}
\end{equation} Further, from \eqref{eq:new4.1},
\eqref{eq:new4.4}, \eqref{eq:new4.10}, and
\eqref{eq:new4.23}, we have
\begin{equation}\label{eq:new4.26}
\Lambda_1=\frac{\tau}{4},
\end{equation}
\begin{equation}\label{eq:new4.27}
\Lambda_2=\frac{1}{48}\Big(\frac{1}{2}\tau^2+3|R|^2\Big),
\end{equation}
\begin{equation}\label{eq:new4.28}
192\Lambda_3=-27|\nabla
R|^2+32\Big(-\frac{7}{2}\hat{R}+\mathring{R}+\frac{\tau^3}{16}+\frac{9}{8}\tau|R|^2\Big).
\end{equation}
From \eqref{eq:new4.26} and \eqref{eq:new4.27}, taking account of
Corollary \ref{cor:4.3}, we see that $\hat{R} - 2\mathring{R}$ is constant. Thus, from Proposition \ref{prop:4.4},
it follows that $\hat{R}$ and $\mathring{R}$ are both constant, and
hence $|\nabla R|^2$ is also constant. Thus, $a$, $b$, and
$c$ are the real roots of the equation
    \begin{equation}\label{new:eq}
    t^3+\frac{\tau}{4}
    t^2+\frac{1}{32}\Big(\frac{5}{6}\tau^2-|R|^2\Big)t-192\Big(\hat{R}-\tau|R|^2-\frac{1}{4}\tau^3\Big)=0
    \end{equation}
    at each point $p\in M$ and hence, $a$, $b$ and $c$ can be expressed in terms of constant-valued functions
     $\tau$, $|R|^2$, $\hat{R}$ and $\mathring{R}$  at each point of $M$, respectively.
 Therefore, $M$ is a
4-dimensional 2-stein curvature homogeneous manifold, and hence $M$
is a locally symmetric manifold by virtue of (\cite{SV}, pp.281).
Further, taking account of the result \cite{Se} and Proposition
\ref{prop4.1}, we may show the following.

\begin{theorem}\label{thm:new4.4}
A 4-dimensional asymptotic harmonic manifold up to order 3 is locally flat or locally isometric to a
rank one symmetric space.
\end{theorem}
Thus, from Theorem \ref{thm:new4.4}, the refinement of the Walker's result follows immediately \cite{W}.

\subsection{5-dimensional asymptotic harmonic manifolds }

First, let $M=(M,g)$ be a 5-dimensional asymptotic harmonic
manifold up to order 2. Then, $M$ is a 2-stein manifold with
constant $|R|^2$. From Corollary \eqref{cor:4.2} and \ref{eq:new4.5}
with $m=5$, we see that $M$ satisfies the equality
\begin{equation}\label{eq:new4.29}
2\mathring{R}_{ij}-\hat{R}_{ij} = \frac{\tau}{100}(9|R|^2 -
\tau^2)g_{ij}.
\end{equation}
Hence, transvecting \eqref{eq:new4.29} with $g_{ij}$, we have
\begin{equation}\label{eq:new4.30}
2\mathring{R}-\hat{R} = \frac{\tau}{20}(9|R|^2-\tau^2).
\end{equation}
Next, let $M = (M,g)$ be a 5-dimensional asymptotic harmonic
manifold up to order 3. Then, from Proposition \ref{prop:4.4} and
\eqref{eq:new4.30}, we see that $\hat{R}$, $\mathring{R}$ and
$|\nabla R|^2$ are constant on $M$ (\cite{W}, Proposition 3.1).
From \eqref{eq:new4.10} with $m=5$, in addition to the
equalities \eqref{eq:new4.29} and \eqref{eq:new4.30}, we have the
following equality
\begin{equation}\label{eq:new4.31}
315\Lambda_3 = -27|\nabla R|^2 - 112\hat{R} + 32\mathring{R} +
\frac{32}{25}\tau^3 + \frac{144\tau}{5}|R|^2.
\end{equation}
Thus, from \eqref{eq:new4.30} and \eqref{eq:new4.31}, we have
\begin{equation}\label{eq:new4.32}
27|\nabla R|^2 + 96\hat{R} - \frac{12}{25} \tau^3 - 36\tau |R|^2 = -
315\Lambda_3.
\end{equation}

Now, we recall the following result of Nikolayevsky (\cite{Ni}, Proposition 1).
\begin{proposition}\label{prop:5.3}
A 5-dimensional 2-stein manifold $M=(M,g)$ is either of constant
sectional curvature or locally homothetic to the symmetric space
$SU(3)/SO(3)$ or to its noncompact dual $SL(3)/SO(3)$.
\end{proposition}

In this section, we give a brief review on Proposition \ref{prop:5.3} under a slightly
more general setting from the view point of Question A. We may note
that the following result (\cite{Ni}, Proposition 4) plays an
essential role in the proof of Proposition \ref{prop:5.3}.
\begin{proposition}\label{prop:5.5}
Let $M=(M,g)$ be a 5-dimensional 2-stein manifold. Then, at each point $p\in M$, there exists an orthonormal basis $\{e_i\}$ such that
\begin{equation*}
\begin{gathered}
R_{1212}=R_{1313}=R_{2323}=R_{2424}=R_{3434}=\mu-\nu,\quad R_{1414}=\mu-4\nu,\\  R_{1515}= R_{4545}=\mu, \quad
R_{2525}=R_{3535}=\mu-3\nu,\\
R_{1234}=\nu,\quad R_{1235}=\sqrt{3}\nu,\quad R_{1324}=-\nu,\quad R_{1325}=\sqrt{3}\nu,\\
R_{1423}=-2\nu,\quad R_{2425}=\sqrt{3}\nu,\quad R_{3435}=-\sqrt{3}\nu
\end{gathered}
\end{equation*}
and all the other components of $R$ vanish up to sign.
\end{proposition}
\noindent From Proposition \ref{prop:5.5}, by direct calculations, we have
\begin{equation}\label{eq:5.31}
\tau=-20\mu+30\nu,
\end{equation}
\begin{equation}\label{eq:new5.33}
R_{iabc}R_{jabc}=(8\mu^2-24\mu\nu+60\nu^2)\delta_{ij},
\end{equation}
and hence,
\begin{equation}\label{eq:5.32}
|R|^2=40\mu^2-120\mu\nu+300\nu^2.
\end{equation}

\noindent {{Further we can obtain the following:
\begin{equation}\label{eq:new5.34}
\check{R}_{ij}=(-32\mu^3+144\mu^2\nu-384\mu\nu^2+360\nu^3)\delta_{ij},
\end{equation}
\begin{equation}\label{eq:5.33}
\hat{R}_{ij}=(16\mu^3-72\mu^2\nu+360\mu\nu^2-600\nu^3)\delta_{ij},
\end{equation}
and hence, $\hat{R}=80\mu^3-360\mu^2\nu+1800\mu\nu^2-3000\nu^3$,
\begin{equation}\label{eq:5.34}
\mathring{R}_{ij}=(12\mu^3-54\mu^2\nu+18\mu\nu^2-30\nu^3)\delta_{ij}
\end{equation}
and hence, $\mathring{R}=60\mu^3-270\mu^2\nu+90\mu\nu^2-150\nu^3$.}}
Thus, from \eqref{eq:5.31} and \eqref{eq:5.32},  we see that $\mu$
and $\nu$ are represented in terms of the constant valued functions
$\tau$ and $|R|^2$ at each point $p\in M$, and hence, $\mu$ and
$\nu$ are constant on $M$. Therefore, $M$ is curvature homogeneous.
From \eqref{eq:new4.14} with $m=5$, taking account of \eqref{eq:5.31}$\sim$\eqref{eq:5.34}, we have
\begin{equation}\label{eq:5.35}
|\nabla R|^2=1680\mu \nu^{2}.
\end{equation}
Thus, from \eqref{eq:5.35}, it follows that $M$ is locally symmetric
if and only if $\mu = 0$ or $\nu = 0$. Here, if $\nu = 0$, then,
from Proposition \ref{prop:5.5}, it follows that $M$ is a space of
constant sectional curvature $-\mu$. Now, we assume that $\nu \neq
0$. Then, by applying the second Bianchi identity to the curvature
form obtained by making use of {Proposition \ref{prop:5.5}}, we
may check that $M$ is locally symmetric (and hence, $\mu=0$), and
further that $M$ is locally homothetic to the symmetric space
$SU(3)/SO(3)$ or to its noncompact dual $SL(3)/SO(3)$ (\cite{Ni},
pp.32$ \sim$ pp.34). Thus, we have Proposition \ref{prop:5.3}.

We now show that any 5-dimensional Riemannian manifold $M=(M,g)$
which is locally homothetic to the symmetric space $SU(3)/SO(3)$
(resp. $SL(3)/SO(3)$) with a fixed canonical Riemannian metric is
never an asymptotic harmonic manifold up to order 3. In order to do
this, without loss of generality, it suffices to establish it in the
case where the Riemannian manifold $M$ is locally homothetic to the
symmetric space $SL(3)/SO(3)$ equipped with the metric given by
(\cite{Ni}, pp.34). Now, we assume that $M$ is an asymptotic
harmonic manifold up to order 3. Then, we may easily check that
$\nu<0$ for $M$. Since $\nabla R=0$ and $\mu=0$ hold on $M$,  from
\eqref{eq:new4.32}, taking account of
\eqref{eq:5.31}$\sim$\eqref{eq:5.34}, we have

\begin{equation}\label{eq:5.36}
\Lambda_{3} = 1984 \nu^{3}.
\end{equation}

On the other hand, choosing an orthonormal basis $\{e_i\} =
\{e_1 = x, e_2,e_3,e_4,e_5\}$ of the tangent space $T_{p}M$ at any
point $p \in M$ satisfying the condition in the Proposition
\ref{prop:5.5} and calculating the equality in the condition
$H'_{3}$ by making using of the orthonormal basis $\{e_{i}\}$, we
have also

\begin{equation}\label{eq:5.37}
\Lambda_{3} = 2012 \nu^{3}.
\end{equation}

Thus, from \eqref{eq:5.36} and \eqref{eq:5.37}, it must follow that $\nu=0$. But,
this is a contradiction. Summing up the above arguments, we have finally the following.
\begin{theorem}\label{thm:5.6}
Let $M=(M,g)$ be a 5-dimensional asymptotic harmonic manifold up to order 3.
Then $M$ is a space of constant sectional curvature.
\end{theorem}
From Theorem \ref{thm:5.6}, we have immediately following
(\cite{Ni}, Theorem 1).
\begin{corollary}\label{corollary:5.7}
A 5-dimensional harmonic manifold is a space of constant sectional curvature.
\end{corollary}
Corollary \ref{corollary:5.7} gives an affirmative answer to the Lichnerowicz conjecture
(refined version by Ledger) for 5-dimensional case.
 \vskip0.2cm

\noindent{\bf Remark 5} The result that the symmetric space
$SU(3)/SO(3)$ (resp. $SL(3)/SO(3)$) is not asymptotic harmonic
manifold up to order 3 can be also obtained by taking account of the
fact that $SU(3)/SO(3)$ (resp. $SU(3)/SL(3)$) is not 3-stein
(\cite{CGW}, pp.58). We here give another explicit proof for the
same result by making use of the curvature identities on
5-dimensional Riemannian manifolds derived from the universal
curvature identity on 6-dimensional Riemannian manifolds.

\section{Concluding remarks}

Based on the discussions in the previous sections, while grappling
with the Lichnerowicz conjecture for 6-dimensional case, it seems
effective to find an orthonormal basis at each point of a
6-dimensional 2-stein manifold such as the Singer-Thorpe basis for
the 4-dimensional case and the Nikolayevsky basis for the
5-dimensional case. As an approach to the Lichnerowicz conjecture
for the 6-dimensional case, it also seems worthwhile to provide the
universal curvature identity on the 8-dimensional Riemannian
manifold through a method similar to the 4- and 5-dimensional cases
and further the curvature identities on the 6- and
7-dimensional Riemannian manifolds derived from the obtained
universal curvature identities.

Lastly, we shall explain a reason why we introduced the notion of
asymptotic harmonic manifolds.
 As mentioned in the beginning of \textsection 5, there are several equivalent definitions for harmonic manifolds.
 One of them is the one defined in terms of the characteristic function $f=f(\Omega)$, where $\Omega=\frac{1}{2}s^2$, $s=d(p,q)$ for $q\in U_p$ ($U_p=U_p(x^1,x^2,\cdots, x^m)$
 denoting a sufficiently small normal coordinate neighborhood centered at each point $p\in M$). The characteristic function plays an important role in the geometry of harmonic manifolds.
 We refer to {\cite{Be, Li44, Ta}} for more details on the characteristic functions. From these observations,
concerning Question A, it is natural to discuss the relationships
between the constants $\{H_n\}_{n=1,2,\cdots}$ and
$\{f^{(n)}(0)\}_{n=1,2,\cdots}$. Here, we denote by `` $\ ' {}$ ''
the derivative with respect to the variable $\Omega$. Now, let
$M=(M,g)$ be an $m$-dimensional harmonic manifold with the
characteristic function $f=f(\Omega)$. Then, it is known that
between the constants $\{\Lambda_1,\Lambda_2,\Lambda_3\}$ and
$\{f'(0),f''(0),f'''(0)\}$, the following relations hold \cite{Li44,
Ta}:
\begin{equation}\label{eq:6.1}
\Lambda_1=-\frac{3}{2}f'(0),\quad \Lambda_2= -\frac{45}{8}f''(0),\quad\Lambda_3=-315f'''(0),
\end{equation}
Lichnerowicz \cite{Li44} has proved the following.
\begin{theorem}\label{thm:6.1}
In any $m$-dimensional harmonic manifold $M=(M,g)$, the characteristic function $f=f(\Omega)$ satisfies the inequality
\begin{equation}\label{eq:6.2}
f'(0)^2+\frac{5}{2}(m-1)f''(0)\leqq0.
\end{equation}
The equality sign is valid if and only if $M$ is of constant sectional curvature.
\end{theorem}
From \eqref{eq:new4.1} and \eqref{eq:new4.4}, taking account of \eqref{eq:6.1} and \eqref{eq:6.2},
we can see that the above Theorem \ref{thm:6.1} is generalized as follows:
\begin{theorem}\label{thm:6.2}
Let $M=(M,g)$ be an $m$-dimensional asymptotic harmonic manifold up
to order 2. Then $M$ satisfies the inequality
\begin{equation}\label{eq:6.3}
\Lambda_1^2-(m-1)\Lambda_2\leqq0
\end{equation}
The equality sign is valid if and only if  $M$ is of constant sectional curvature $\frac{\tau}{m(m-1)}$.
\end{theorem}
Tachibana \cite{Ta} has proved the following
\begin{theorem}\label{thm:6.3}
Any $2n$-dimensional harmonic K\"ahler manifold $M=(M,J,g)$ satisfies the inequality
\begin{equation}\label{eq:6.4}
f'(0)^2+\frac{5(n+1)^2}{n+7}f''(0)\leqq0.
\end{equation}
and the equality sign is valid if and only if $M$ is of constant holomorphic sectional curvature.
\end{theorem}
From \eqref{eq:new4.1} and \eqref{eq:new4.4}, taking account of \eqref{eq:6.1} and \eqref{eq:6.2},
we can see that the above Theorem \ref{thm:6.3} is generalized as follows:
\begin{theorem}\label{thm:6.4}
Let $M=(M,J,g)$ be a $2n$-dimensional asymptotic harmonic K\"ahler
manifold up to order 2. Then $M$ satisfies the inequality
\begin{equation}\label{eq:6.5}
\Lambda_1^2-\frac{2(n+1)^{2}}{n+7}\Lambda_2\leqq0
\end{equation}
and the equality sign is valid if and only if $M$ is of constant
holomorphic sectional curvature $\frac{\tau}{n(n+1)}$.
\end{theorem}

Similarly, from \eqref{eq:new4.1}, \eqref{eq:new4.4} and
\eqref{eq:new4.10}, taking account of \eqref{eq:6.1}, we can see
that the corresponding generalizations for the results (\cite{Wa1},
Theorem 5.2) and (\cite{Wa2}, Theorem 5.5) are
obtained.

Taking account of the discussions in
the present paper and \cite{CGW}, concerning the Question A, we
obtain that if the dimension is 4 then the least integer of series
is not greater than 3 and if the dimension is 5 then the least
integer of series is 3. Based on the arguments developed the the
following question will naturally arise:
 \vskip0.2cm

\noindent{\bf Question B.} For any integer $m(m \geqq 6)$,
does there exist the least integer $K(m)$ such that any
$m$-dimensional asymptotic harmonic manifold up to order $k(k \geqq
K(m))$ is necessary harmonic ?

\section*{Acknowledgements}
\noindent This research was supported by Basic Science Research
Program through the National Research Foundation of Korea(NRF)
funded by the Ministry of Education(2014053413).

\end{document}